\documentclass[11pt]{article}
\usepackage{amsfonts}
\usepackage{color}
\usepackage{graphics}

\usepackage{cite}
\usepackage{latexsym}
\usepackage{amsmath}
\usepackage{amssymb}
\usepackage[dvips]{epsfig}
\usepackage{amscd}

 \usepackage{color}
\newtheorem{theorem}{Theorem}[section]
\newtheorem{remark}{Remark}[section]

\newtheorem{lemma}[theorem]{Lemma}

\newcommand{\ti}{\tilde}

\def\pf{{\it Proof.}  }

\newcommand{\pa}{\partial}

\newcommand{\thatsall}{\hfill$\Box$}
\newcommand{\bi}{\bibitem}

\newcommand{\dis}{\displaystyle}

\newcommand{\bt}{\begin{theorem}}
\newcommand{\bl}{\begin{lemma}}
\newcommand{\el}{\end{lemma}}
\newcommand{\et}{\end{theorem}}
\newcommand{\ga}{\gamma}

\renewcommand{\b}{\beta  }

\newcommand{\te}{\theta}

\newcommand{\al}{\alpha}
\newcommand{\de}{\delta}

\newcommand{\la}{\label}
\newcommand{\si}{\sigma}
\newcommand{\ka}{\kappa}

\newcommand{\bn}{\begin{eqnarray}}
\newcommand{\en}{\end{eqnarray}}
\newcommand{\bnn}{\begin{eqnarray*}}
\newcommand{\enn}{\end{eqnarray*}}

\newcommand{\bnnn}{\begin{eqnarray*}}
\newcommand{\ennn}{\end{eqnarray*}}
\newcommand{\ben}{\begin{enumerate}}
\newcommand{\een}{\end{enumerate}}

\newcommand{\ba}{\begin{aligned}}
\newcommand{\ea}{\end{aligned}}
\newcommand{\be}{\begin{equation}}
\newcommand{\ee}{\end{equation}}

\newcommand\abs[1]{\left\lvert #1 \right\rvert}
\newcommand\norm[1]{\left\lVert #1 \right\rVert}

\def\xix{\int_{\O} }

\def\xiT{\int_0^T}

\def\O{\Omega}

\def\r{\mathbb{R}}

\makeatletter
\@addtoreset{equation}{section}
\makeatother

\makeatletter
\@addtoreset{equation}{section}
\makeatother

\title{ Global Existence  of
 Strong Solutions to  Compressible Navier-Stokes System  with
Degenerate Heat Conductivity in Unbounded Domains\thanks{Partially supported by Undergraduate Research Fund of BNU 2017-150,
  Email addresses: kexinli98@163.com (K. Li),   shuxl03@hotmail.com (X. Shu), xjxu@bnu.edu.cn(X. Xu)}}

\date{}
\author{   Kexin L{\small I}, Xuanlin S{\small HU},  Xiaojing X{\small U}
 \\[3mm] {\normalsize
 School of Mathematical Sciences,}\\
 {\normalsize Beijing Normal University,} \\
  {\normalsize Beijing 100875,  P. R. China }  }

\date{ }

\begin{document}
\maketitle

\begin{abstract} In one-dimensional unbounded domains, we  prove global existence of strong solutions     to the compressible Navier-Stokes system for a viscous and heat conducting ideal polytropic gas,
when the viscosity is constant and the heat conductivity $\ka$
depends on the temperature $\te$ according to $\ka=\bar\ka \te^\b (\b>0)$.
Note that the conditions imposed on  the  initial data are the same as those of the constant  heat conductivity case ([Kazhikhov, A. V.   Siberian Math.
J. 23 (1982), 44-49]) and can be arbitrarily large.
Therefore, our result generalizes  Kazhikhov's result  for the constant heat conductivity case to the degenerate and nonlinear one. \end{abstract}


Math Subject Classification: 35Q35; 76N10.

\section{Introduction}
The compressible Navier-Stokes system  describing the  one-dimensional motion of a
viscous heat-conducting    gas
can be written in the Lagrange variables in the following form (see \cite{ba,se})\be \la{1.1}
v_t=u_{x},
\ee
\be \la{1.2}
u_{t}+P_{x}=\left(\mu\frac{u_{x}}{v}\right)_{x},\ee
\be  \la{1.'3} \left(e+\frac{u^2}{2}\right)_{t}+ (P
u)_{x}=\left(\kappa\frac{\theta_{x}}{v}+\mu\frac{uu_x}{v}\right)_{x}
 ,
\ee
where   $t>0$ is time,  $x\in \O\subset \r=(-\infty,+\infty)$ denotes the
Lagrange mass coordinate, and  the unknown functions $v>0,u,$ $\theta>0,e>0,$ and $P$ are,  respectively, the specific volume of the gas,  fluid velocity,   absolute temperature,   internal energy, and  pressure. In general, $P,e,$  the viscosity $\mu$, and heat conductivity $\ka$ are functions of $\te$ and
$v.$  In this paper, we
consider ideal polytropic gas, that is, $P$ and $e$ satisfy \be  \la{1.3''}  P =R \theta/{v},\quad e=c_v\theta +\mbox{ const.},
\ee
where $R$ (specific gas constant) and $c_v $ (heat
capacity at constant volume) are both  positive constants. For $\mu$ and $\ka,$
we  consider   the case where $\mu$ and $\ka$ are proportional to (possibly different) powers of  $\te:$  \be   \la{1.3'}\mu=\ti\mu \te^\ga, \quad \ka=\ti\ka \te^\beta, \ee
where $\ti\mu,\ti\ka>0$ and $\ga,\beta\ge 0$ are constants.

The system \eqref{1.1}-\eqref{1.3'} is imposed on
  the initial condition
\be \la{1.4} (v(x,0),u(x,0),\theta(x,0))=(v_0(x),u_0(x),\theta_0(x)),\quad x\in
\O,\ee
and   three types of far-field and  boundary conditions:

1) Cauchy problem
  \be \la{1.5}  \O=\r,\, \lim\limits_{|x|\rightarrow
  \infty}(v(x,t),u(x,t),\theta(x,t))=(1,0,1),\quad t>0;\ee

   2)   boundary and far-field conditions for $ \O=(0,\infty),$
\be \la{1.6}  u(0,t)=0,\, \theta_x(0,t)=0,\,\lim\limits_{x\rightarrow
\infty}(v(x,t),u(x,t),\theta(x,t))=(1,0,1),\quad t>0; \ee

 3)   boundary and far-field conditions for $ \O=(0,\infty),$
\be \la{1.7} u(0,t)=0,\, \theta (0,t)=1,\, \lim\limits_{x\rightarrow
\infty}(v(x,t),u(x,t),\theta(x,t))=(1,0,1), \quad t>0. \ee

According to the Chapman-Enskog expansion for the first level of approximation in kinetic theory, the viscosity $\mu$
and heat conductivity $\ka$ are functions of temperature alone (\cite{cip,cc}). In particular, if the intermolecular
potential varies as $r^{-a},$  with intermolecular distance  $r$, then $\mu$ and $\ka$ are both
proportional to the power $(a+4)/(2a)$ of the temperature, that is, \eqref{1.3'} holds with $\ga=\beta=(a+4)/(2a).$ Indeed, for Maxwellian molecules $(a=4),$ the dependence is linear, while for elastic spheres $(a\rightarrow\infty),$ the dependence is like $\te^{1/2}.$

For constant coefficients $(\ga=\b =0),$ Kazhikhov and Shelukhin
\cite{ks} first obtained    the global existence
     of solutions  in bounded domains  for large initial data.  From then on,
significant progress has been made on the mathematical aspect of the initial and
initial
boundary value problems, see \cite{az1,az2,az3,akm,ji4,ji3,kaw} and the references therein. For the Cauchy problem  \eqref{1.1}-\eqref{1.5} and the initial boundary value
problems
  \eqref{1.1}-\eqref{1.4}   \eqref{1.6}   and \eqref{1.1}-\eqref{1.4}   \eqref{1.7}
  (in unbounded domains), Kazhikhov \cite{ka1} (also cf. \cite{akm,ji4}) obtains the global existence of strong solutions.

Much effort has been made to generalize this approach to other cases and in
particular to models satisfying \eqref{1.3'}. However,  there are few results for the case that $\ga=\beta$,  partially because of the   difficulty introduced in such relations, which
lead to possible degeneracy and strong nonlinearity in viscosity and heat diffusion. As
a first step in this direction,  for bounded domain $\O,$ Jenssen-Karper \cite{jk1}  proved the global existence of   weak solutions to \eqref{1.1}--\eqref{1.3'}
under the assumption that $\ga=0$ and $\beta\in(0,3/2).$ Later, for $\ga=0$ and $\beta\in(0,\infty),$  Pan-Zhang \cite{pz1} obtain the global strong solutions. However, their methods depend heavily on the boundedness of the domain $\O.$ In this paper, we will  prove the global existence of strong solutions to  the Cauchy problem  \eqref{1.1}-\eqref{1.5} and the initial boundary value
problems
  \eqref{1.1}-\eqref{1.4}   \eqref{1.6}   and \eqref{1.1}-\eqref{1.4}   \eqref{1.7}.  That is,
  our main result is as follows.
 \begin{theorem}\la{thm11}  Suppose that \be \ga=0,  \quad \beta> 0,\ee
 and  that the initial data $(v_0,u_0,\theta_0)$ satisfy
\be\la{c5} v_0-1,u_0,\theta_0-1\in H^1(\O),\,\, \inf_{x\in \O}v_0(x)>0,\, \,
\inf_{x\in \O}\theta_0(x)>0 ,\ee
and   are compatible with \eqref{1.6},  \eqref{1.7}.
Then there exists a unique global
strong solution
  $(v ,u ,\theta )$   with positive  $v(x, t)$ and  $\theta(x, t)$ to the initial-boundary-value problem \eqref{1.1}-\eqref{1.5},  or \eqref{1.1}-\eqref{1.4}  \eqref{1.6}, or
  \eqref{1.1}-\eqref{1.4}  \eqref{1.7}   satisfying for any $T>0,$ \be
  \la{3k}\begin{cases}   v-1, \,u,\, \theta-1 \in L^\infty(0,T;H^1(\O) ), \\ v_t\in
  L^\infty(0,T;L^2 (\O))\cap L^2(0,T;H^1(\O)), \\ u_t,\,\theta_t,\,v_{xt},\,u_{xx},\,\theta_{xx} \,\in
  L^2(\O\times(0,T)).\end{cases}\ee
  \end{theorem}
\begin{remark} Our Theorem \ref{thm11} can be regarded as a natural  generalization   of  Kazhikhov's result (\cite{ka1}) where he considered the constant viscosity case ($\ga=\b=0$) to the degenerate and nonlinear one that $\ga=0,\b>0.$  \end{remark}

\begin{remark} Although for $\ga=\b=0,$ the large-time behavior of the strong solution to the initial-boundary-value problem \eqref{1.1}-\eqref{1.5},  or \eqref{1.1}-\eqref{1.4}  \eqref{1.6}, or
  \eqref{1.1}-\eqref{1.4}  \eqref{1.7} has been proved just recently (see \cite{ji1,ji2,ll}), their methods cannot be applied directly to our case  due to the the   degeneracy and nonlinearity of the heat conductivity since $\b>0.$ It is interesting to study the large-time behavior of the strong solutions for $\b>0.$ This will be left for future. \end{remark}

We now comment on the analysis of this paper. Compared with  the constant viscosity case ($\ga=\b=0$) (\cite{ka1}), the main difficulty comes from the   degeneracy and nonlinearity of the heat conductivity due to the fact that $\b>0.$ The key observations are as follows: First, we modify slightly the idea due to Kazhikhov \cite{ka1} to prove the lower and upper bounds of $v.$ Next, multiplying the energy equation \eqref{1.'3} by $\te^{-2}(\te^{-1}-2)_+^p$ and noticing that the domains $\{x\in\O|\te(x,t)<1/2\}$ remain  bounded for all $t\in [0,T]$ (see \eqref{aq1s}), we find that the temperature is indeed bounded from below (see Lemma \ref{lkq4}) which lays a firm foundation for our
further analysis. Finally, to obtain the higher order estimates, we will modify some ideas due to Li-Liang \cite{ll} to obtain the estimates on the $L^2(\O\times(0,T))$-norms of both $u_t$ and $u_{xx}$ (see Lemma \ref{lwq2}) which are crucial for further estimates on the upper bound of the temperature.The whole procedure will be carried out in the next section.

\section{ Proof of Theorem \ref{thm11}}

We first state  the local existence lemma which can be proved by using the principle of compressed mappings (c.f. \cite{kan,nash,tan}).

\begin{lemma} Under the conditions of  Theorem \ref{thm11},  there exists some $T>0$ such that the initial-boundary-value problem \eqref{1.1}-\eqref{1.5},  or \eqref{1.1}-\eqref{1.4}  \eqref{1.6}, or
  \eqref{1.1}-\eqref{1.4}  \eqref{1.7} has a unique
strong solution
  $(v ,u ,\theta )$   with positive  $v(x, t)$ and  $\theta(x, t)$   satisfying \eqref{3k}.\end{lemma}

Then, to finish the proof of   Theorem \ref{thm11}, it only remains to obtain  a priori   estimates (see \eqref{cqq1},  \eqref{aq3s}, \eqref{vvx1}, \eqref{qu3g}, \eqref{tq7} below) the constants in
which depend only on $T$ and the data of the problem.
The estimates make it possible to continue the local solution  to the whole interval $[0,T]$.

Next, without loss of generality, we assume that $\ti\mu= \ti\ka=R=c_v=1 .$
We modify slightly the idea due to Kazhikhov \cite{ka1} to prove the lower and upper bounds of $v.$

\begin{lemma} \la{lwq3} There exists a positive constant $C$ such that
	\be\ba\la{cqq1}
	C^{-1}\le v(x,t)\le C,
	\ea\ee where (and in what follows)   $C $  denotes a
	generic positive constant
	depending only on $T, \b,\|(v_0-1,u_0,\theta_0-1)\|_{H^1(\O)},
	\inf\limits_{x\in \O}v_0(x),$ and $ \inf\limits_{x\in \O}\theta_0(x).$\end{lemma}

\pf
First,
using \eqref{1.1}, \eqref{1.2}, and \eqref{1.3''}, we rewrite
\eqref{1.'3} as
\be\la{tee1}\theta_{t}+  \frac\theta v
u_{x}= \left(\frac{\theta^\b\theta_{x}}{v}\right)_{x}+ \frac{u_{x}^{2}}{v}.\ee Multiplying  (\ref{1.1}) by
 $ 1- {v}^{-1} $, (\ref{1.2})  by $
u$,    and  (\ref{tee1}) by
$  1- {\theta}^{-1} $,  and adding them altogether, we obtain
\bnn  \ba&(u^2/2+ (v-\log v-1)+ (\theta-\log
\theta-1))_t+ \frac{u^2_x}{v\theta}
+ \frac{\theta_x^2}{v\theta^2}\\&=\left(\frac{  u u_x}{v}-\frac{  u
\theta}{v}\right)_x + u_x+ \left((1-\theta^{-1})
\frac{\theta_x}{v}\right)_x,\ea\enn  which together with \eqref{1.5} or \eqref{1.6}
or \eqref{1.7} yields
  \be\ba\la{etvt} \sup_{0\le t\le T}\xix \left(  \frac{u^2}{2}+ (v-\log
v-1)+ (\theta-\log \theta-1)\right)dx+\int_0^T W(s)ds \le e_0,  \ea\ee
where
\be\ba W(t)\triangleq\xix \left(\frac{\theta^\b\theta_{x}^2}{v\te^2}+\frac{u_{x}^{2}}{v\te}
\right)(x,t)dx, \ea\ee and \bnn e_0\triangleq\xix\left(\frac{u_0^2}{2}+ (v_0-\log v_0-1)+ (\theta_0-\log
\theta_0-1)\right)dx.\enn

Next, for any $x\in \O, $ denoting $ N=[x],$   we have by \eqref{etvt}
\be\la{ja1}\int_{N}^{N+1}(v-\log v-1)+(\theta-\log\theta-1)dx\le e_0,\ee
which together with Jensen's inequality yields \be\la{ab} \al_1 \le \int_N^{N+1}v(x,t)dx\le \al_2, \al_1\le \int_N^{N+1}\te(x,t)dx\le \al_2,\ee   where $0<\al_1<\al_2$  are two roots of  \bnn x-\log x-1=e_0.\enn
Moreover, it follows from \eqref{ja1} that for any $t>0,$ there exists some
 $b_N(t)\in[N,N+1]$ such that\bnn (v-\log v-1+\theta-\log\theta-1)(b_N(t),t)\le e_0,\enn
which implies  \be\la{ab'} \al_1\le v(b_N(t),t)\le \al_2, \,\al_1\le \te(b_N(t),t)\le \al_2. \ee

Letting $\si\triangleq\dfrac{u_x}{v}-\dfrac{\te}{v}=(\log v)_t-\dfrac{\te}{v}$, we   write \eqref{1.2} as
\be\la{1.2'} u_t=\si_x.\ee
Integrating  \eqref{1.2'} over $[N,x]\times[0,t]$ leads to
\bnn  \int_N^x(u(y,t)- u_0)dy = \log v-\log v_0-\int_0^t\frac{\te}{v}ds-\int_0^t\si(N,s)ds,\enn
which gives
\be\la{DY}\ba v(x,t)  &= D_N(x,t) Y_N(t) \exp\left\{\int_0^t\frac{\te}{v}ds\right\}, \ea\ee
where
\bnn D_N(x,t)\triangleq v_0(x)\exp\left\{\int_{N}^x \left(u(y,t)-u_0(y)\right)dy\right\},
\enn and
\bnn Y_N(t)\triangleq \exp\left\{\int_0^t\si (N,s)ds\right\}.\enn
Thus, it follows from  \eqref{DY} that
\be\la{DY'} v(x,t)=  D_N(x,t) Y_N(t) \left(1+ \int_0^t\frac{\te(x,\tau)}{ D_N(x,\tau) Y_N(\tau)}d\tau \right).\ee
Since $$\left| \int_N^x  \left(u(y,t)-u_0(y)\right)dy\right|\leq \left(\dis\int_N^{N+1} u^2dy\right)^{1/2}+\left(\dis\int_N^{N+1} u_0^2dy\right)^{1/2}\leq C,$$  we have \be\la{D}C^{-1}\leq D_N(x,t)\leq C,\ee where and in what follows,  $C$ is a constant independent of $N$.

Moreover, integrating \eqref{DY'} with respect to $x$ over $[N,N+1]$ gives
\bnn   \frac{1}{ Y_N(t)}\int_N^{N+1}v(x,t)dx= \int_N^{N+1} D_N(x,t) \left(1+ \int_0^t\frac{\te(x,\tau)}{ D_N(x,\tau) Y_N(\tau)}d\tau \right)dx,\enn which yields
\be\la{DD1}\ba C^{-1}\le  \frac{1}{ Y_N(t)}\le C+C\int_0^t \frac{1}{ Y_N(\tau)}d\tau,\ea\ee where we have used \eqref{ab}, \eqref{D},  and the following simple fact:
\bnn \int_N^{N+1} \frac{\te(x,\tau) D_N(x,t)}{ D_N(x,\tau)}dx\le C\int_N^{N+1} \te dx\le C.\enn
Applying  Gr\"onwall's inequality to \eqref{DD1}  gives
\bnn\la{Y} C^{-1}\le\frac{1}{ Y_N(t)}\le C,\enn which together with \eqref{DY'} and  \eqref{D} implies that for $(x,t)\in [N,N+1]\times [0,T],$
 \be\la{eq1}    C^{-1}\le v(x,t)\le C +C \int_0^t \max\limits_{x\in [N,N+1]}\te(x,t) dt.\ee

Then, direct computation gives
\bnn\ba&\left|\te^{\frac{\beta+1}{2}}(x,t)- \te^{\frac{\beta+1}{2}}(b_N(t),t) \right| \\ &=\left|\int_{b_N(t)}^x\left(\te^{\frac{\beta+1}{2}}\right)_xdx\right|\\& =\frac{\beta+1}{2}\left|\int_{b_N(t)}^x\frac{\te^{\beta/2}\te_x}{\sqrt{\te}}dx\right| \\&\le \frac{\beta+1}{2}\left(\int_N^{N+1} \frac{ \te^\beta \te_x^2}{\te^2 v} dx\right)^{1/2}\left(\int_N^{N+1}  {\te v} dx\right)^{1/2}\\ &\le C W^{1/2}(t)\max\limits_{x\in [N,N+1]}v^{1/2}(x,t)  , \ea\enn which together with \eqref{ab'}  yields that for any $t>0,$
\bnn
\max\limits_{x\in [N,N+1]}\te(x,t)\le C\left(1+W (t)\max\limits_{x\in [N,N+1]}v(x,t)\right)   .
\enn Putting this into (\ref{eq1}) and
using    Gr\"onwall's inequality yields that
\bnn \max\limits_{(x,t)\in [N,N+1]\times[0,T]}v(x,t)\le C,\enn
which combined with \eqref{eq1}    and the fact that $C$ is independent of $N$ gives \eqref{cqq1}.    The proof of Lemma \ref{lwq3} is finished. \thatsall

Now we are in a position  to obtain the lower bound of the temperature $\te.$
\begin{lemma} \la{lkq4} There exists a positive constant $C$ such that for all $(x,t)\in \O\times[0,T],$
	\be  \la{aq3s}  \te(x,t) \ge C^{-1}. \ee \end{lemma}

\pf First, denoting by $$(\te>2)(t) =\{x\in\O\mid\te(x,t)>2\},$$ and  $$(\te<1/2)(t) =\{x\in\O\mid\te(x,t)<1/2\},$$
we get by \eqref{etvt}
\bnn\ba e_0&\ge    \int_{(\te<1/2)(t)} (\te-\log \te-1)dx+ \int_{(\te>2)(t)} (\te-\log \te-1)dx \\& \ge (\log 2-1/2)\abs{(\te<1/2)(t)}+(1-\log 2 )\abs{(\te>2)(t)}
\\& \ge (\log 2-1/2)\left(\abs{(\te<1/2)(t)}+\abs{(\te>2)(t)}\right),\ea\enn
which shows that for any $t\in [0,T],$
\be\la{aq1s}\abs{(\te<1/2)(t)}+\abs{(\te>2)(t)}\le \frac{2e_0}{2\log 2-1}.\ee

Next, for $p>2$, multiplying \eqref{tee1} by $\theta^{-2}(\theta^{-1}-2)_+^p$ with $(\theta^{-1}-2)_+\triangleq\max\{\theta^{-1}-2,0\}  $ and integrating over $\O$, we obtain
 \bnn\ba &\dfrac{1}{p+1}\left(\int_\O (\theta^{-1}-2)_+^{p+1}dx\right)_t+\int_\O \dfrac{u_x^2}{v\theta^2}(\theta^{-1}-2)_+^pdx\\
	&\le\int_\O \dfrac{u_x}{v\theta}(\theta^{-1}-2)_+^pdx\\
	&\le\dfrac{1}{2}\int_\O \dfrac{u_x^2}{v\theta^2}(\theta^{-1}-2)_+^pdx+2\int_\O \dfrac{1}{v}(\theta^{-1}-2)_+^pdx \\
	&\le\dfrac{1}{2}\int_\O \dfrac{u_x^2}{v\theta^2}(\theta^{-1}-2)_+^pdx+C\left(\int_\O (\theta^{-1}-2)_+^{p+1}dx\right)^{\frac{p}{p+1}},
\ea\enn
where in the last inequality we have used \eqref{aq1s}. Thus,  we have \bnn \|(\theta^{-1}-2)_+\|_{L^{p+1}(\O)}^{p}\left(\|(\theta^{-1}-2)_+\|_{L^{p+1}(\O)}\right)_t\le C\|(\theta^{-1}-2)_+\|_{L^{p+1}(\O)}^{p}, \enn with $C$ independent of $p$. This in particular implies that there exists some positive constant $C$ independent of $p$ such that
\be \la{aq2s}
\sup_{0\le t\le T}\norm{(\theta^{-1}-2)_+}_{L^{p+1}(\O)}\le C.\ee
 Using \eqref{aq1s}, letting $p\rightarrow \infty$ in \eqref{aq2s}  shows \bnn
\sup_{0\le t\le T}\norm{(\theta^{-1}-2)_+}_{L^\infty(\O)}\le C,\enn which proves \eqref{aq3s} and finishes the proof of Lemma \ref{lkq4}. \thatsall

For further uses, we need the following estimates on the $L^2(\O\times (0,T))$-norm of $u_x.$

\begin{lemma} \la{lkq1} There exists a positive constant $C$ such that
\be   \la{lmm3}\xiT\xix \left(u_{x}^2+\te^{-1}\te_x^2\right)dxdt   \le C.\ee
\end{lemma}

\pf First, by  \eqref{etvt}, we have
\be \la{ng2'} \ba  \int_0^T \max_{x\in \O}\te^{ 1+\beta } dt&\le C \int_0^T \max_{x\in \O}(\te-2)_+^{1+\beta}dt+C\\&\le C\int_0^T\int_{\Omega}(\te-2)_+^{\beta}\left|\te_x\right|dxdt+C
\\&\le C\xiT \int_{(\te>2)(t)} \te ^{\beta}|\te_x|dx dt+C\\&\le C\xiT \left(\int_{(\te>2)(t)} \te^{\beta+2} dx\right)^{1/2}\left(\int_\Omega\te^{\beta-2}\te_x^2dx \right)^{1/2} dt+C\\&\le
 C\xiT \left(\max_{x\in \O}\te^{ 1+\beta }\right)^{1/2}\left(\int_\Omega\te^{\beta-2}\te_x^2dx \right)^{1/2} dt+C
\\&\le \frac12 \xiT \max_{x\in \O}\te^{ 1+\beta }dt+C \xiT \int_\Omega\te^{\beta-2}\te_x^2dx  dt+C,\ea\ee where in the fourth  inequality we have used
    \be \la{nep1}\sup\limits_{0\le t
\le T}\int_{(\te>2)(t)}\theta dx \le C \sup\limits_{0\le t
\le T}\int_{\O}(\theta-\log \theta-1)dx\le C , \ee due to \eqref{etvt}.
Combining  \eqref{ng2'}, \eqref{etvt}, and \eqref{cqq1} yields
\be  \ba  \int_0^T \max_{x\in \O}\te^{ 1+\beta } dt&\le C,
\ea\ee which implies
\be \la{ng2} \ba \int_0^T \max_{x\in \O}\te dt\le C \int_0^T \max_{x\in \O}\left(1+\te^{ 1+\beta }\right) dt&\le C.
\ea\ee

Next, integrating the momentum equation  \eqref{1.2}  multiplied by $u$ with respect to $x$  over $\O,$ after  integrating by parts, we obtain
\bnn \ba & \frac12\left(\xix u^2dx\right)_t+\xix \frac{u_x^2}{v}dx \\&=\xix \frac{\te }{v}u_xdx \\&=\xix \frac{\left(\te-1\right)}{v}u_xdx-\xix \frac{ \left(v-1\right)u_x}{v}dx   \\&\le C\xix \left(\te-1\right)^2dx+C\xix  \left(v-1\right)^2dx +\frac12\xix \frac{u_x^2}{v}dx \\&\le C\int_{(\te>2)} \te^2dx+C+C\xix  \left(v-1\right)^2dx +\frac12\xix \frac{u_x^2}{v}dx   \\&\le C\max_{x\in\O}\te  +C +\frac12\xix \frac{u_x^2}{v}dx   , \ea\enn where in the last inequality we have used \eqref{nep1}, \eqref{etvt}, and \eqref{cqq1}.
Combining this  with \eqref{ng2} gives \be \la{bg1}\ba  \xiT\xix u_x^2dxdt\le C.\ea\ee

Finally, if $\b\ge 1,$ we have \be  \la{qs1}  \xiT\xix \te^{-1}\te_x^2dxdt   \le C\xiT\xix \te^{\b-2}\te_x^2dxdt\le C. \ee
If $\beta\in(0,1),$  for $0<2\al<1,$ multiplying \eqref{tee1}   by $\te^{\al-1}(\te^\al-2^\al)_+$    and integration by parts gives
\be\nonumber\ba &  (1-2\al)\int_{(\te>2)(t)}  \frac{\te^{\b+2\al-2}\te_x^2}{v }dx  +\xix \frac{u_x^2}{v }\te^{\al-1}(\te^\al-2^\al)_+dx\\& =   \frac{1}{2\al}\left(\xix  (\te^\al-2^\al)_+^2dx\right)_t+2^\al (1-\al )\int_{(\te>2)(t)} \frac{\te^{\b+\al-2}\te_x^2}{v }dx\\&\quad+\xix   \frac{  \te^\al (\te^\al-2^\al)_+u_x}{v}dx \\ & \le   \frac{1}{2\al}\left(\xix  (\te^\al-2^\al)_+^2dx\right)_t+ \frac{1-2\al}{2}\int_{(\te>2)(t)}  \frac{\te^{\b+2\al-2}\te_x^2}{v }dx\\&\quad +C ( \al )\int_{(\te>2)(t)}\frac{\te^{\b -2}\te_x^2}{v }dx+C\xix   u_x^2dx+C+C\max_{x\in\O}\te\int_{(\te>2)(t)} \te dx,\ea\ee
which together with \eqref{bg1}, \eqref{ng2},  \eqref{nep1},  \eqref{aq1s}, and Gronwall's inequality yields
\bnn \xiT\xix   \te^{\b+2\al-2}\te_x^2 dxdt\le C(\al). \enn In particular, combining this where we  choose $2\al=1-\b  \in(0,1),$   \eqref{qs1}, and  \eqref{bg1}  proves \eqref{lmm3}. The proof of Lemma \ref{lkq1} is completed. \thatsall

\begin{lemma} \la{lwq1} There exists a positive constant $C$ such that

 \be\la{vvx1} \ba
\sup_{0\le t \le T}\xix  v_x^2  dx+\int_0^T\xix  v_x^2\te dxdt\le C.\ea\ee

\end{lemma}

\pf
Rewriting momentum equation \eqref{1.2} as
 \be\ba\la{mom1}
\left(\frac{v_x}{v}\right)_t=u_t+\left(\frac{\te}{v}\right)_x,
\ea\ee due to \bnn\ba
\left(\frac{v_t}{v}\right)_x=\left(\frac{v_x}{v}\right)_t ,
\ea\enn
we multiply  \eqref{mom1}  by $\dfrac{v_x}{v}$  to get
 \be\ba\la{mom2}
\frac{1}{2}\left[\left(\frac{v_x}{v}\right)^2\right]_t&=\frac{v_x}{v}u_t+\frac{v_x}{v}\left(\frac{\te}{v}\right)_x
\\&=\left(\frac{v_x}{v}u\right)_t-u(\log v)_{xt}+\frac{v_x\te_x}{v^2}-\frac{v_x^2\te}{v^3}
\\&=\left(\frac{v_x}{v}u\right)_t-\left[u(\log v)_{t}\right]_x+\frac{u_x^2}{v}+\frac{v_x\te_x}{v^2}-\frac{v_x^2\te}{v^3}.\ea\ee
Integrating \eqref{mom2} over $\O\times[0,T]$ and using \eqref{lmm3},  one has
 \bnn\ba\la{mom3}
&\sup_{0\le t \le T}\xix\left[\frac{1}{2}\left(\frac{v_x}{v}\right)^2-\frac{v_x}{v}u\right]dx
+\int_0^T\xix\frac{v_x^2\te}{v^3}dxdt\\
&\le C+\int_0^T\xix  \frac{v_x\te_x}{v^2}dxdt\\
&\le C+\frac12\int_0^T\xix \frac{v_x^2\te}{v^3}dxdt+ C\int_0^T\xix   \te^{-1} \te_x^2dxdt\\
&\le C+\frac12\int_0^T\xix \frac{v_x^2\te}{v^3}dxdt.
\ea\enn   This in particular implies \eqref{vvx1} due to the following simply fact: \bnn\la{mom4}\ba
\xix \frac{v_x}{v}udx\le \frac{1}{4}\xix \left(\frac{v_x}{v}\right)^2dx+ C.
\ea\enn
Thus, the proof of Lemma \ref{lwq1} is finished.\thatsall
\begin{lemma} \la{lwq2}There exists a positive constant $C$ such that

 \be \la{qu3g} \ba \sup_{0\le t\le T}\xix  u_x^2  dx+\xiT\xix (u_{xx}^2+  u_t^2)dxdt \le C.  \ea\ee

\end{lemma}

\pf
  First, we rewrite the momentum equation  \eqref{1.2}  as
\be\ba\la{mom5} u_t-\frac{u_{xx}}{v}=-\frac{u_xv_x}{v^2}-\frac{\te_x}{v } +\frac{\te v_x}{v^2}.\ea\ee
Multiplying both sides of \eqref{mom5} by $u_{xx},$ and integrating the resultant equality in $x$ over $\O ,$ one has
\be\la{qu2}\ba
&\frac{1}{2}\frac{d}{dt}\xix u_x^2dx+\xix\frac{u_{xx}^2}{v} dx\\
&\le \left|\xix\frac{u_xv_x}{v^2}u_{xx}dx\right|+\left|\xix\frac{\te_x}{v }u_{xx}dx\right| +\left|\xix\frac{\te v_x}{v^2}u_{xx}dx\right|\\
&\le \frac{1}{4}\xix\frac{u_{xx}^2}{v} dx+C\xix\left(u_x^2v_x^2+v_x^2\te^2 +\te_x^2\right)dx.
\ea\ee
Direct computation yields that for any $\de>0,$
\be\la{qu1}\ba &\xix \left(u_x^2v_x^2+v_x^2\te^2+\te_x^2\right)dx\\ &\le C\left(\max_{x\in\O } u_x^2+\max_{x\in\O } \te^2\right)\xix  v_x^2dx +\xix \te_x^2dx\\ &\le C\max_{x\in\O } u_x^2+C\max_{x\in\O } (\te-2)_+^2  +C +\xix \te_x^2dx\\ &\le \delta \xix u_{xx}^2dx+ C(\delta) \xix u_x^2dx  +C +C\xix \te_x^2dx,\ea\ee
where in the last inequality  we have used
\be\ba\la{ux2}
\max_{x\in\O }u_x^2&\le \xix  \left|\left(u_x^2\right)_x\right|dx\\&\le 2\left(\xix  u_{xx}^2d x\right)^{1/2}\left(\xix  u_x^2dx\right)^{1/2} \\&\le \delta \xix  u_{xx}^2dx+ C(\delta) \xix  u_x^2dx,
\ea\ee and \bnn\ba \max_{x\in\O } (\te-2)_+^2&=\max_{x\in\O }\left(\int_x^\infty \pa_y(\te-2)_+(y,t)  dy\right)^2\\&\le \left(\int_{(\te>2)(t)}  |\te_y| dy\right)^2\\&\le C\xix \te_x^2dx.\ea\enn
Putting \eqref{qu1} into \eqref{qu2} and choosing $\de$ suitably small yields
\be\la{qu2'}\ba
& \frac{d}{dt}\xix u_x^2dx+\xix\frac{u_{xx}^2}{v} dx\\
&\le   C   \xix u_x^2dx  +C +C_1\xix \te^\b\te_x^2dx.
\ea\ee

Next, motivated by \cite{ll}, we integrate  \eqref{tee1} multiplied  by $(\theta-2)_{+}\triangleq \max\{\theta-2,0\}$ over $\O$ to get
    \begin{equation}\begin{split} \label{1.1-1}  &\frac{1}{2} \left(\int_{\O}
    (\theta-2)_{+}^{2}dx\right)_t+ \int_{(\te>2)(t)}\frac{ \te^\b\theta_x^{2}}{v}dx
    \\  & =  - \int_{\O}\frac{\theta}{v}u_{x}(\theta-2)_{+} dx+ \int_{\O}
    \frac{u_{x}^{2}}{v}(\theta-2)_{+}dx\\  & \le   C \max_{x\in \O}\te\left( \int_{\O} (\theta-2)^2_{+} dx+   \int_{\O}
     {u_{x}^{2}}  dx\right). \end{split}\end{equation}

Noticing that \bnn\ba  \int_{\O}  \te^\b \te_x^2dx&\le   \int_{(\te>2)(t)}  \te^\b \te_x^2dx+  \int_{ (\te\le 2)(t)}  \te^\b \te_x^2dx \\ &\le C  \int_{(\te>2)(t)}\frac{ \te^\b\te_x^2}{v}dx+ C \int_{(\te\le 2)(t)} \frac{ \te^{\b-2}\te_x^2}{v}dx, \ea\enn
we deduce from \eqref{1.1-1} that
    \bnn\ba   &\frac{1}{2} \left(\int_{\O}
    (\theta-2)_{+}^{2}dx\right)_t+ C_2\int_{\Omega }    \te^\b \theta_x^{2} dx
    \\    & \le C+ C \xix  \frac{ \te^{\b-2}\te_x^2}{v}dx + C \max_{x\in \O}\te \int_{\O} (\theta-2)^2_{+} dx+ C \max_{x\in \O}\te \int_{\O}
     {u_{x}^{2}}  dx. \ea \enn
    Adding this multiplied by $2C_2^{-1}C_1 $ to \eqref{qu2'} together with Gronwall's inequality gives
 \be \la{qu3} \ba \sup_{0\le t\le T}\xix (u_x^2+(\te-2)_+^2) dx+\xiT\xix (u_{xx}^2+  \te^\b \te_x^2)dxdt \le C,  \ea\ee
 which together with \eqref{mom5}, \eqref{qu2}, and \eqref{qu1} yields \eqref{qu3g} and finishes the proof of Lemma \ref{lwq2}. \thatsall

\begin{lemma}\la{tq8}There exists a positive constant $C$ such that \be\la{tq7}\ba \sup_{0\le t\le T}\xix  \te_x^2dx+\int_0^T\xix \left( \te_t^2+\te_{xx}^2\right)dxdt\le C.\ea\ee\end{lemma}
\pf
  First, multiplying \eqref{tee1} by
$  \te^{\b}\te_t$ and integrating the resultant equality over $\O$ yields
\bnn\ba & \xix  \te^\b\te_t^2dx+\xix \frac{ \te^{\b+1}\theta_tu_x}{v  }dx\\&=\xix \te^{\b}\te_t\left(\frac{\te^\b\theta_{x}}{v}\right)_{x}dx+\xix \frac{ \te^{\b }\theta_tu_x^2}{v }dx  \\& =-\xix \frac{\te^\b\theta_{x}}{v}\left(\te^{\b}\te_t\right)_{x}dx+\xix \frac{ \te^{\b }\theta_tu_x^2}{v }dx   \\& =-\xix \frac{\te^\b\theta_{x}}{v}\left(\te^{\b}\te_x\right)_{t}dx+\xix \frac{ \te^{\b }\theta_tu_x^2}{v }dx  \\&=-\frac{1}{2}
\xix \frac{\left((\te^\b\theta_{x})^2\right)_t}{v}dx+\xix \frac{ \te^{\b }\theta_tu_x^2}{v }dx \\&=-\frac{1}{2} \left(\xix \frac{(\te^\b\theta_{x})^2}{v}dx\right)_t
-\frac{1}{2}\xix \frac{(\te^\b\theta_{x})^2u_x}{v^2}dx+\xix \frac{ \te^{\b }\theta_tu_x^2}{v }dx ,\ea\enn
which gives
\be\ba\la{lm8eq2} & \xix  \te^\b\te_t^2dx+ \frac{1}{2} \left(\xix \frac{(\te^\b\theta_{x})^2}{v}dx\right)_t\\&=- \frac{1}{2}\xix \frac{(\te^\b\theta_{x})^2u_x}{v^2}dx-\xix \frac{ \te^{\b+1}\theta_t u_x}{v  }dx+\xix \frac{ \te^{\b }\theta_tu_x^2}{v }dx\\&\le C\max_{x\in\O} (|u_x|\te^{\b/2})\xix \left(\te^{3\b/4}\te_x\right)^2dx +\frac{1}{2}\xix  \te^\b\te_t^2dx \\&\quad+C\xix  \te^{\b+2}u_x^2dx+C\xix \te^\b u_x^4dx\\&\le C \xix  \te^{2\b}\te_x^2dx\xix  \te^\b\te_x^2dx  +C \max_{x\in\O} (\te^{\b+2}+ \te^\b u_x^2)+\frac{1}{2}\xix  \te^\b\te_t^2dx+C  \\&\le  C \xix  \te^{2\b}\te_x^2dx\xix  \te^\b\te_x^2dx+\frac{1}{2}\xix  \te^\b\te_t^2dx+C \max_{x\in\O}(  \te^{2\b+2}+  u_x^4)+C  \ea\ee
due to \eqref{qu3}.

Next, it follows from  \eqref{ux2}  and \eqref{qu3}    that
\bnn\la{tq3} \int_0^T \max_{x\in\O}  u_x^4dt\le C,\enn which together with \eqref{lm8eq2},  the Gronwall inequality, and \eqref{qu3} leads to
\be\ba\la{tebtex}  \sup_{0 \le t\le T}\xix  \left(\te^\b\theta_{x}\right)^2 dx+\int_0^T\xix  \te^\b\te_t^2dxdt\le C, \ea\ee where we have used \be \la{teb1}\max_{x\in\O} \te^{2\beta+2}\le C+C\xix (\te^\b\te_x)^2dx.\ee
Combining \eqref{teb1} with \eqref{tebtex}  implies that for all $(x,t)\in \O\times [0,T],$
\be\la{teup} \te(x,t)\le C.\ee

Then, combining \eqref{aq3s} and \eqref{tebtex} leads to
\be\ba\la{lm9eq2}  \sup_{0 \le t\le T}\xix  \theta_{x}^2 dx+\int_0^T\xix  \te_t^2dxdt\le C. \ea\ee

Finally, it follows from \eqref{tee1} that
\bnn\ba \frac{\te^\b\te_{xx}+\b\te^{\b-1}\te_x^2 }{v}=\frac{\left(\te^\b\te_x\right)_x}{v}=  \frac{ \te^\b\te_xv_x}{v^2}- \frac{u_x^2}{v}+ \frac{ \te u_x}{v}+\te_t,\ea\enn
which together with \eqref{aq3s}, \eqref{teup},  \eqref{lm9eq2}, \eqref{qu3},   and \eqref{ux2} gives
\bnn\la{tq4}\ba  \int_0^T\xix \te_{xx}^2dxdt  & \le C\int_0^T\xix \left(\te_x^4+\te_x^2v_x^2+u_x^4+u_x^2+\te_t^2\right)dxdt\\ & \le C+ C\int_0^T\max_{x\in\O} \left(\te_x^2 +u_x^2 \right) dt\\ & \le C+ C\int_0^T   \xix \te_x^2 dx   dt+ \frac12\int_0^T  \xix \te_{xx}^2 dx   dt.\ea\enn

 Combining this with \eqref{lm9eq2}
gives \eqref{tq7} and finished the proof of Lemma \ref{tq8}.\thatsall

 \end{document}